\documentclass[francais]{smfart}

\usepackage[dvips]{epsfig}

\usepackage[francais]{babel}

\usepackage{smfthm}

\def\R{\ifmmode{{\rm I}\hskip -3pt {\rm R}}
    \else{\hbox{${\rm I}\hskip -3pt {\rm R}$}}\fi} 
\def\P{\ifmmode{{\rm I}\hskip -3pt {\rm P}}
    \else{\hbox{${\rm I}\hskip -3pt {\rm P}$}}\fi} 
\def\N{\ifmmode{{\rm I}\hskip -3pt {\rm N}}
    \else{\hbox{${\rm I}\hskip -3pt {\rm N}$}}\fi} 
\def\Z{\ifmmode{{\rm Z}\hskip -4.8pt {\rm Z}}
    \else{\hbox{${\rm Z}\hskip -4.8pt {\rm Z}$}}\fi} 
\def\C{\ifmmode{{\rm C}\hskip-4.8pt\vrule height5.8pt\hskip6.3pt}
    \else{\hbox{${\rm C}\hskip-4.8pt\vrule height5.8pt\hskip6.3pt$}}\fi}
\def\Q{\ifmmode{{\rm Q}\hskip-5.0pt\vrule height6.0pt depth 0pt\hskip6pt}
    \else{\hbox{${\rm Q}\hskip-5.0pt\vrule height6.0pt depth 0pt\hskip6pt$}}\fi}
\def\D{\ifmmode{{\rm I}\hskip -3pt {\rm D}}
    \else{\hbox{${\rm I}\hskip -3pt {\rm D}$}}\fi} 

\newtheorem{theorem}{Th\'eor\`eme}[section]
\newtheorem{lemma}[theorem]{Lemme}
\newtheorem{proposition}[theorem]{Proposition}
\newtheorem{corollary}[theorem]{Corollaire}
\newtheorem{rappel}[theorem]{}

\theoremstyle{definition}
\newtheorem{definition}[theorem]{D\'efinition}
\newtheorem{remark}[theorem]{Remarque}

\newtheorem{notation}[theorem]{Notation}
\newtheorem{question}[theorem]{Question}
\newtheorem{example}[theorem]{Exemple}
\newtheorem{conjecture}[theorem]{Conjecture}
\newtheorem{exercice}[theorem]{Exercice}

\newcommand{\e}{{\rm e}^}

\newcommand{\cal}{\mathcal}

\newcommand{\sn}{{\rm sn}\,}
\newcommand{\cn}{{\rm cn}\,}
\newcommand{\dn}{{\rm dn}\,}
\newcommand{\cP}{\check{\P}\C^2}

\newcommand{\bi}{\begin{itemize}}
\newcommand{\ei}{\end{itemize}}
\newcommand{\be}{\begin{enumerate}}
\newcommand{\ee}{\end{enumerate}}

\newcommand{\bpf}{\begin{proof}}
\newcommand{\epf}{\end{proof}}

\newcommand{\bt}{\begin{theorem}}
\newcommand{\et}{\end{theorem}}
\newcommand{\brap}{\begin{rappel}}
\newcommand{\erap}{\end{rappel}}
\newcommand{\bnt}{\begin{notation}}
\newcommand{\ent}{\end{notation}}
\newcommand{\bd}{\begin{definition}}
\newcommand{\ed}{\end{definition}}
\newcommand{\ble}{\begin{lemma}}
\newcommand{\ele}{\end{lemma}}
\newcommand{\bpr}{\begin{proposition}}
\newcommand{\epr}{\end{proposition}}
\newcommand{\bre}{\begin{remark}}
\newcommand{\ere}{\end{remark}}
\newcommand{\bco}{\begin{corollary}}
\newcommand{\eco}{\end{corollary}}
\newcommand{\beq}{\begin{equation}}
\newcommand{\eeq}{\end{equation}}
\newcommand{\bq}{\begin{question}}
\newcommand{\eq}{\end{question}}
\newcommand{\beqn}{\begin{eqnarray*}}
\newcommand{\eeqn}{\end{eqnarray*}}
\newcommand{\bex}{\begin{example}}
\newcommand{\eex}{\end{example}}
\newcommand{\ber}{\begin{exercice}}
\newcommand{\eer}{\end{exercice}}
\newcommand{\sct}{\section}

\newcommand{\sk}{\smallskip}
\newcommand{\bk}{\bigskip}
\newcommand{\nk}{\noindent}
\newcommand{\pl}{\partial}

\newcommand{\fr}{\frac}
\newcommand{\bcj}{\begin{conjecture}}
\newcommand{\ecj}{\end{conjecture}}

\author{Luc Pirio et Jean-Marie Tr\'epreau}

\address{UMR 7586,  175 rue du Chevaleret, 75013 Paris}

\email{pirio@math.jussieu.fr ; trepreau@math.jussieu.fr}

\title{Une famille de 5-tissus plans exceptionnels}

\begin{document}

\begin{abstract}
Le $5$-tissu de Bol, d\'ecouvert en 1936, est rest\'e longtemps le seul 
exemple av\'er\'e de tissu plan exceptionnel. En 2002, Robert \cite{Ro}
et Pirio \cite{Pi} ont trouv\'e un $9$- et des $6$- et $7$-tissus 
plans exceptionnels,  en relation avec l'\'equation fonctionnelle
de Spence-Kummer du trilogarithme. Plus r\'ecemment encore, Pirio \cite{Pi2} a 
d\'ecouvert trois nouveaux $5$-tissus exceptionnels,
tr\`es simples. 
Dans cet article, nous pr\'esentons une famille
\`a un param\`etre de $5$-tissus plans exceptionnels, qui contient 
ces derniers comme cas limites. 

Dans leur pr\'esentation la plus naturelle, les tissus de cette 
famille sont compos\'es d'un syst\`eme harmonique
de quatre faisceaux de droites parall\`eles, d\'efinis 
respectivemant par les \'equations  $\{x={\rm cte}\}$,
$\{y={\rm cte}\}$, $\{x+y={\rm cte}\}$ et $\{x-y={\rm cte}\}$,
et d'un feuilletage dont les feuilles sont les courbes de niveau 
de la fonction ${\rm sn}_k x \, {\rm sn}_k y$, o\`u ${\rm sn}_k$ est une 
fonction elliptique de Jacobi. 

\end{abstract}

\maketitle

\sct{Introduction}

Pour fixer les id\'ees, nous nous pla\c{c}ons dans une situation locale, dans
la cat\'egorie analytique. On note $(x,y)$ le point courant de $\C^2$ et on choisit 
un point de base $(x_0,y_0)$ ; fonction voudra dire germe de fonction analytique en $(x_0,y_0)$, etc...

\'Etant donn\'e une fonction $u$ telle que $du(x_0,y_0)\neq 0$, 
on note $\cal{F}(u)$ le feuilletage dont les feuilles sont les courbes 
de niveau $\{u(x,y)= {\rm cte}\}$ de la fonction $u$. On dit que $u$
est une fonction de d\'efinition du feuilletage. Si $f$ est une fonction
d\'efinie au voisinage de $u(x_0,y_0)$ et si $f'(u(x_0,y_0))\neq 0$,
la fonction $f(u):=f\circ u$ d\'efinit le m\^eme feuilletage que $u$.

Soit $d\geq 3$ un entier. Un $d$-tissu en $(x_0,y_0)$ est une famille 
non ordonn\'ee de $d$ feuilletages $\cal{F}(u_k)$, deux \`a deux transverses, {\em i.e.}
tels que :
$$
\text{si $j,k = 1,\ldots, d\,$ et $j\neq k$}, \qquad du_j(x_0,y_0)\wedge du_k(x_0,y_0) \neq 0.
$$
On notera ce tissu :
\beq
\label{tissu}
\cal{T}(u_1,\ldots,u_d) := \{\cal{F}(u_1),\ldots,\cal{F}(u_d)\}.
\eeq
On dira que deux (germes de) $d$-tissus $\cal{T}_0$  en $(x_0,y_0)$ et 
$\cal{T}_1$ en $(x_1,y_1)$ sont {\em \'equivalents} s'il existe un 
(germe de) diff\'eomorphisme $\phi$, $\phi(x_0,y_0)=(x_1,y_1)$, 
qui envoie $\cal{T}_0$ sur $\cal{T}_1$. 

\sk
La th\'eorie des tissus est n\'ee \`a la fin des ann\'ees 1920. 
Le livre \cite{BB} de Blaschke et Bol
dresse le bilan des r\'esultats fondamentaux obtenus 
entre 1927 et 1938. Les r\'esultats que nous 
rappelons maintenant sont tous pr\'esent\'es dans \cite{BB},
sauf la d\'ecouverte r\'ecente de nouveaux tissus plans exceptionnels
et certains travaux oubli\'es des ann\'ees 1937--1939 que nous 
\'evoquons \`a la fin de cette introduction.

Apr\`es la guerre et surtout apr\`es 1970, 
de nouveaux r\'esultats ont \'et\'e 
obtenus, en particulier en dimension ou en codimension 
plus grandes. Citons seulement, \`a titre d'exemples, 
les articles de Chern et Griffiths \cite{CG}, 
d'Atkivis et Goldberg  \cite{AG} et de H\'enaut \cite{He}, 
dont on pourra consulter les bibliographies. 

\sk
Une {\em relation ab\'elienne} du tissu $\cal{T}(u_1,\ldots,u_d)$
est une relation fonctionnelle de la forme :
$$
\sum_{j=1}^d \phi_j(u_j)\, du_j = 0,
$$
o\`u $\phi_j$ est une fonction d\'efinie au voisinage de $u_j(x_0,y_0)$. 
On identifie une telle relation au $d$-uplet $(\phi_1(u_1),\ldots,\phi_d(u_d))$. 
Les relations ab\'eliennes 
forment un espace vectoriel, ici sur le corps des complexes.
La dimension $\rho$ de cet espace est appel\'ee {\em le rang du tissu}.
Elle  v\'erifie {\em l'in\'egalit\'e de Bol} :
\beq
\label{borne}
\rho \leq \rho(d), \qquad \rho(d):= \fr{1}{2}(d-1)(d-2).
\eeq
On a par exemple $\rho(3)=1$, $\rho(4)=3$, $\rho(5)=6$, $\rho(6) = 10$...
\bd
Un $d$-tissu de rang $\rho(d)=(d-1)(d-2)/2$ est dit de rang maximal.
\ed
On \'ecrira plut\^ot une  relation ab\'elienne sous la forme int\'egr\'ee :
$$
\sum_{j=1}^d f_j(u_j)= {\rm cte},
$$
o\`u $f_j$ est une fonction d\'efinie au voisinage de $u_j(x_0,y_0)$. 
En faisant cela, on introduit des relations triviales, les 
constantes $f=(f_1,\ldots,f_d)\in \C^d$ ; le rang du tissu est 
la dimension de l'espace des $d$-uplets $(f'_1(u_1),\ldots,f'_d(u_d))$.

\sk
Par exemple, consid\'erons $d\geq 3$ formes lin\'eaires $\,l_1,\ldots,l_d$
sur $\C^2$, deux \`a deux lin\'eairement ind\'ependantes. 
Le tissu $\cal{T}(l_1,\ldots,l_d)$, compos\'e de $d$ faisceaux 
diff\'erents de droites parall\`eles, est de rang maximal.
On peut le voir en remarquant que, pour tout $n \in \N^*$,
le syst\`eme $(l_1^n,\ldots,l_d^n)$ est de rang $\min (n+1,d)$ dans 
l'espace des polyn\^omes homog\`enes de degr\'e $n$ en $(x,y)$,
qui est un espace de dimension $n+1$. Si $1 \leq n \leq d-2$, il 
existe donc $d-n-1$ relations lin\'eaires ind\'ependantes entre 
$l_1^n,\ldots,l_d^n$. En sommant par rapport \`a $n$,
on obtient $(d-2) + (d-3) + \cdots + 1 = (d-1)(d-2)/2$
relations ab\'eliennes.

C'est d\'ej\`a un exercice 
plus difficile de montrer directement qu'un tissu compos\'e
de $d$ faisceaux de droites non n\'ecessairement parall\`eles
est de rang maximal. C'est un cas tr\`es particulier d'un r\'esultat 
que nous rappelons maintenant.

\sk
On note $\cP$ l'espace projectif 
dual de $\P\C^2$ : un point $p\in \cP$
repr\'esente une droite $l(p)$ de $\C\P^2$ ;
\`a tout point $z\in \P\C^2$ est associ\'e le faisceau de droites 
de sommet $z$, donc une droite $l_z$ de $\cP$.

Soit $C \subset \cP$ une courbe alg\'ebrique r\'eduite de degr\'e 
$d$. Si $z$ est un \og point g\'en\'erique \fg\, de $\P\C^2$,
la droite $l_z$ de $\cP$ coupe la courbe $C$ en 
$d$ points distincts $p_1,\ldots,p_d$ de $\cP$.
En associant \`a $z$ les $d$ droites $l(p_1),\ldots,l(p_d)$
(elles passent par $z$), on d\'efinit un \og $d$-tissu avec singularit\'es \fg\, 
sur $\C\P^2$. Au point g\'en\'erique de $\C\P^2$, on obtient 
un $d$-tissu au sens usuel. C'est un {\em tissu en droites}, c'est-\`a-dire 
que ses feuilles sont des morceaux de droites. 

En g\'en\'eral, le tissu obtenu s'interpr\`ete comme 
le tissu des tangentes \`a la courbe duale de $C$.
Dans le cas particulier o\`u la courbe $C$ est la r\'eunion de $d$ droites 
distinctes, on obtient un tissu de $d$ faisceaux
de droites. 

\bd
Un tissu obtenu de la fa\c{c}on qu'on vient de d\'ecrire 
est dit alg\'ebrique.
\ed
Les deux r\'esultats suivants sont fondamentaux. Le premier est une interpr\'etation,
due \`a Blaschke, du th\'eor\`eme d'addition d'Abel en termes de tissus :
\bt[\cite{BB}]
Tout tissu alg\'ebrique est de rang maximal.
\et
Le second, d\^u \`a Blaschke et Howe, est souvent cit\'e 
comme \og le th\'eor\`eme d'Abel inverse \fg\, :
\bt[\cite{BB}]
\label{abel-inv}
Tout tissu en droites de rang maximal est un tissu alg\'ebrique.
\et 
\bd
Un tissu \'equivalent \`a un tissu alg\'ebrique est dit 
alg\'ebrisable. Un tissu de rang maximal non alg\'ebrisable est
dit exceptionnel.
\ed
Par exemple et c'est trivial, tout $3$-tissu $\cal{T}(u_1,u_2,u_3)$
de rang $\rho(3)=1$ est alg\'ebrisable. En effet, si 
$f_1(u_1) + f_2(u_2) + f_3(u_3) = {\rm cte}$ est une
relation ab\'elienne non triviale du tissu, on v\'erifie 
facilement que $f'_j(u_j(x_0,y_0))\neq 0$, $\, j=1,2,3$. On peut donc 
prendre  $f_1(u_1), f_2(u_2)$ et $-f_3(u_3)$
comme fonctions de d\'efinition des feuilletages. En changeant de notation, 
la relation ab\'elienne s'\'ecrit $u_3 = u_1 + u_2$ et,
en prenant $(u_1,u_2)$ comme coordonn\'ees locales, on voit 
que le tissu est \'equivalent au tissu $\cal{T}(x,y,x+y)$.

\sk
De m\^eme et ce n'est pas trivial, {\em tout $4$-tissu de rang maximal est alg\'ebrisable}. 
Sous une forme diff\'erente, ce r\'esultat est essentiellement d\^u \`a Lie.

\sk
Le premier exemple de tissu exceptionnel fut donn\'e par Bol  
en 1935~\cite{BB} ; c'est le $5$-tissu suivant :
\beq
\label{Bol}
\cal{B} = \cal{T}( \, x, \, y, \, \fr{y}{x}, \, \fr{1-y}{1-x}, \, \fr{x - xy}{y - xy} \,).
\eeq
Il est compos\'e de quatre faisceaux de droites dont les sommets 
sont en position g\'en\'erale dans $\P\C^2$
(on ne peut pas en extraire  trois points align\'es) ;
le dernier feuilletage est le faisceau des coniques 
qui passent par les quatre sommets. L'une des
relations ab\'eliennes du tissu de Bol fait intervenir
le dilogarithme.

\sk
L'exemple de Bol est rest\'e longtemps isol\'e.
En 2002, Robert \cite{Ro} et le premier auteur \cite{Pi}
ont montr\'e ind\'ependamment que le $9$-tissu associ\'e \`a 
l'\'equation fonctionnelle de Spence-Kummer du trilogarithme
est exceptionnel, ainsi que certains des $6$- et des $7$-tissus 
qu'on peut en extraire.

\sk
En 2003, le premier auteur \cite{Pi2} a donn\'e les trois 
exemples suivants, tr\`es simples, de $5$-tissus exceptionnels :
\begin{align*}
& \cal{T}( \, x, \, y, \, x+y, \, x-y, \, x^2 + y^2 \, ),     \\
& \cal{T}( \, x, \, y, \, x+y, \, x-y, \, x^2 - y^2 \, ),     \\
& \cal{T}( \, x, \, y, \, x+y, \, x-y, \, \e{x}  + \e{y} \, ).
\end{align*}
Dans cet article, nous pr\'esentons une famille \`a un param\`etre 
de $5$-tissus exceptionnels ainsi que cinq tissus exceptionnels 
\og isol\'es \fg\, 
qui apparaissent comme des tissus limites de la famille.
On retrouve en particulier les trois exemples pr\'ec\'edents.

\sk
Ces tissus sont pr\'esent\'es, ainsi que leurs relations ab\'eliennes,
dans la section \S 3, dont la lecture est suffisante si l'on ne s'int\'eresse
qu'aux r\'esultats. Tous ces tissus sont de la forme :
\beq
\label{Tu}
\cal{T}[u]:= \cal{T}( \, x, \, y, \, x+y, \, x-y, \, u(x,y) \, ).
\eeq
Dans tout l'article, on notera :
\beq
\label{T0}
\cal{T}_0:= \cal{T}( \, x, \, y, \, x+y, \, x-y \, ).
\eeq
Dans la section \S 2, on fait quelques remarques sur les sym\'etries
du $4$-tissu $\cal{T}_0$ et les transformations qui envoient 
un tissu de la forme (\ref{Tu}) sur un tissu de la m\^eme  forme.
Dans la section \S 4, on \'ecrit le syst\`eme d'\'equations
que doit v\'erifier la fonction $u(x,y)$ pour que le tissu
(\ref{Tu}) soit de rang maximal. On le r\'esout sous l'hypoth\`ese 
suppl\'ementaire que le $3$-tissu extrait $\cal{T}(\, x, \, y, \, u(x,y) \,)$
est aussi de rang maximal. On v\'erifie que les seules solutions 
sont, \`a \'equivalence de tissus pr\`es, les tissus 
pr\'esent\'es dans la section~\S 3. 
La section \S 4 est d'une certaine fa\c{c}on inutile, 
mais elle indique en partie la voie qu'on a suivie pour d\'ecouvrir
ces nouveaux tissus exceptionnels.

\sk
En fait, les tissus dont nous allons parler sont d\'ej\`a apparus
dans la litt\'erature, sans que leur caract\`ere exceptionnel 
ait \'et\'e reconnu.

\sk
Avant que Bol obtienne son contre-exemple, Blaschke~\cite{Bl} 
a cherch\'e \`a montrer que tout $5$-tissu de
rang $6$ est alg\'ebrisable. La m\'ethode qu'il propose
g\'en\'eralise une m\'ethode de Poincar\'e, 
qui permet de d\'emontrer assez facilement le r\'esultat de
Lie, que tout $4$-tissu de rang $3$ est alg\'ebrisable.
Cet article de Blaschke et le contre-exemple de Bol 
ont suscit\'e plusieurs travaux, de Bompiani et Bortolotti \cite{Bo}
et de Terracini \cite{Te} en 1937, de Buzano \cite{Bu} \cite{Bu2} en 1939.

Ces auteurs s'int\'eressent \`a la question des $5$-tissus exceptionnels
en relation avec un probl\`eme de g\'eom\'etrie diff\'erentielle 
projective.
Il s'agit de trouver les surfaces (S) de l'espace 
$\P\R^5$ qui admettent cinq syst\`emes de courbes 
($=$ un $5$-tissu) tels que (S) soit tangente 
le long de chaque courbe \`a un hyperplan de $\P\R^5$.
Ils trouvent plusieurs solutions 
\`a ce probl\`eme. {\em On sait maintenant que 
toutes ces solutions correspondent \`a des tissus 
exceptionnels \cite{Pi3}.}

L'article \cite{Bu2} est \'ecrit en termes de tissus
et \og tous nos tissus \fg\, apparaissent dans cet article.
Buzano distingue ces tissus pour certaines 
propri\'et\'es de leurs sous-tissus, mais il
ne voit pas qu'ils sont de rang maximal ; 
essentiellement, il manque \`a chaque fois {\em une relation}
pour que leur nature exceptionnelle soit r\'ev\'el\'ee.
Bien que n'ayons eu connaissance de \cite{Bu2}
qu'apr\`es avoir obtenu nos r\'esultats, il y a 
co\"{\i}ncidence entre certains calculs de 
la section \S 4 et ceux de \cite{Bu2}.

Quoiqu'il en soit, Blaschke et Bol connaissaient 
ces travaux. L'appendice \`a la section \S 30 de \cite{BB}, page 261,
contient un commentaire sur les articles \cite{Bo} et \cite{Te},
dont voici un extrait : 

{\em \og ... ist noch nicht festgestellt,
ob die hierzu geh\"origen Kurvengewebe wirklich den Rang
sechs haben, das w\"are aber leicht nachzupr\"ufen \fg}.

Nous ne savons pas si cette v\'erification facile a jamais \'et\'e 
faite. Un peu plus tard, Bol \'ecrit les r\'esum\'es des 
deux articles de Buzano pour le {\em Zentralblatt}
sans commentaire qui aille au-del\`a de ce que l'auteur 
dit d\'emontrer. En 1985, Chern \cite{Ch} pose comme un 
probl\`eme important et non r\'esolu de savoir si le 
tissu de Bol est le seul $5$-tissu non alg\'ebrisable 
de rang 6. En 2004, apr\`es la d\'ecouverte des tissus 
exceptionnels li\'es au trilogarithme, Griffiths \cite{Gr} 
repose la question, sous une forme modifi\'ee :
les tissus exceptionnels ont-ils toujours \`a voir 
avec les polylogarithmes ? Tous nos exemples 
montrent que la r\'eponse est non.

\sct{Remarques sur les tissus $\cal{T}(\, x,\, y, \, x+y,\, x-y, u(x,y) \,)$}

Le tissu
$$
\cal{T}_0:= \cal{T}(\, x, \, y, \, x+y, \, x-y \,)
$$
est compos\'e de quatre faisceaux de droites parall\`eles,
de sommets les points $\,p_1=[0 : 1 : 0]$, $\,p_2=[1 : 0 : 0]$,
$\, p_3=[1 : -1 : 0]$
et $\, p_4=[1 : 1 : 0]$ de $\P\C^2$. 
En particulier, les sommets sont align\'es, alors que 
les sommets du $4$-tissu en droites extrait du tissu 
de Bol (\ref{Bol}) sont en position g\'en\'erale.

Le tissu $\cal{T}_0$ est encore plus particulier que cela.
Rappelons 
que le birapport d'un syst\`eme $(q_1,q_2,q_3,q_4)$ de 
quatre points distincts d'une droite projective est 
donn\'e par :
$$
b(q_1,q_2,q_3,q_4) 
= 
\fr{t_1 - t_3}{t_1 - t_4} :  \fr{t_2 - t_3}{t_2 - t_4},
$$
o\`u $t_k=\xi_k/\eta_k\in \C\cup\{\infty\}$ et les
$[\xi_k:\eta_k]$ sont les coordonn\'ees homog\`enes des $q_k$
dans un rep\`ere projectif de la droite ; c'est un invariant 
projectif du syst\`eme. 
On v\'erifie que $b(p_1,p_2,p_3,p_4)=-1$. {\em Autrement dit, les sommets 
$p_1,p_2,p_3,p_4$  forment une division harmonique de la droite \`a l'infini.}

\sk
Le tissu $\cal{T}_0$ est \'evidemment invariant par les dilatations
$$
(x,y) \mapsto (kx+x',ky+y'), \qquad (k\in \C^*, \; (x',y')\in \C^2).
$$
Si l'on se restreint \`a $\R^2\subset \C^2$, les feuilles de $\cal{T}_0$
issue de $O=(0,0)$ rencontrent le cercle de centre $O$ et de rayon $1$
aux sommets d'un octogone r\'egulier. Il est clair que le tissu $\cal{T}_0$
est invariant par le  groupe des isom\'etries de cet octogone. 
C'est un groupe \`a 16 \'el\'ements, isomorphe au groupe di\'edral $\D_8$
et qu'on notera abusivement $\D_8$. Il permet de substituer \`a $(x,y)$ un des couples 
suivants :
$$
(\pm x, \pm y), \;\; (\pm y, \pm x),
\;\;
(\pm \fr{x+y}{\sqrt{2}},\pm \fr{x-y}{\sqrt{2}}), \;\; (\pm \fr{x-y}{\sqrt{2}},\pm \fr{x+y}{\sqrt{2}}),
$$
avec deux choix ind\'ependants des signes pour chaque couple.
Il est engendr\'e par les sym\'etries orthogonales 
\beq
\label{sym}
\rho(x,y)=(y,x), \qquad \sigma(x,y)=((x+y)/\sqrt{2},(x-y)/\sqrt{2}),
\eeq
par rapport aux droites d'angles polaires $\pi/4$ et $\pi/8$,
respectivement : $\tau \circ \sigma$ engendre le groupe des huit rotations de 
$\D_8$ et, en conjuguant $\tau$ et $\sigma$ par ces rotations, 
on obtient les huit r\'eflexions de $\D_8$. 

Les lemmes qui suivent nous seront utiles pour reconna\^{\i}tre les tissus 
exceptionnels de la forme (\ref{Tu}) et les classer \`a \'equivalence pr\`es.
\ble
\label{LA}
Soit $\phi$ un diff\'eomorphisme local de $\C^2$. 
Si $\phi$ transforme localement les quatre faisceaux de $\cal{T}_0$ en faisceaux de 
droites, $\phi$ est le germe d'une tranformation projective.
Si $\phi$ conserve localement trois des faisceaux
de $\cal{T}_0$, $\phi$ est le germe d'une dilatation.
\ele
\bpf
Pour la premi\`ere propri\'et\'e, en composant $\phi$ avec une 
transformation projective, on peut supposer que :
$$
\phi(\cal{F}(x)) = \cal{F}(x), \qquad \phi(\cal{F}(y))= \cal{F}(y),
$$
et que
$$
\phi(\cal{F}(x+y))=\cal{F}(x+y) \;\; \text{ ou } \;\; \phi(\cal{F}(x+y))=\cal{F}(x/y),
$$
selon que les sommets des faisceaux $\phi(\cal{F}(x))$, $\phi(\cal{F}(y))$
et $\phi(\cal{F}(x+y))$ sont ou ne sont pas align\'es.
De la premi\`ere propri\'et\'e, il d\'ecoule que $\phi^{-1}$ est de la forme :
$$
\phi^{-1}(x,y) = (X(x),Y(y)).
$$
Si $\phi(\cal{F}(x+y)=\cal{F}(x+y)$, on a :
$$
(\pl/\pl x - \pl/\pl y)(X(x)+Y(y)) = X'(x) - Y'(y) = 0,
$$
d'o\`u $X(x)=kx+x'$, $Y(y)=ky+y'$ avec $k\in \C^*$ : $\phi$ est une dilatation
(\c{c}a d\'emontre la deuxi\`eme partie de l'\'enonc\'e).

Si $\phi(\cal{F}(x+y))=\cal{F}(x/y))$, on a :
$$
( x \pl/\pl x + y\pl/\pl y) (X(x)+Y(y)) = xX'(x) + yY'(y) = 0,
$$
d'o\`u $X(x)= k\log x + x'$, $Y(y)= -k\log y + y'$.
On obtient alors 
$$
\phi(\cal{F}(x-y)) = \cal{F}(X(x) - Y(y)) = \cal{F}(xy),
$$
qui n'est pas un faisceau de droites ; c'est une contradiction.
\epf
\ble
\label{LB}
Le groupe ${\rm Sym}\,(\cal{T}_0)$ des transformations projectives qui 
conservent le tissu $\cal{T}_0$ est engendr\'e par les dilatations de $\C^2$
et le groupe $\D_8$. De plus, si $\phi$ est un diff\'eomorphisme local 
qui conserve le tissu $\cal{T}_0$, $\phi$ est le germe 
d'un \'el\'ement de ${\rm Sym}\,(\cal{T}_0)$.
\ele
\bpf
Une transformation projective $\phi$ qui conserve le tissu $\cal{T}_0$
conserve la famille $\{p_1,p_2,p_3,p_4\}$ des sommets 
des faisceaux de droites de $\cal{T}_0$.
Elle conserve donc la droite \`a l'infini (c'est une tranformation affine)
et induit une permutation $s$ de $1,2,3,4$ : $\phi(p_k)=p_{s(k)}$.
De plus, elle conserve le birapport.
Comme les points $p_1,p_2,p_3,p_4$ forment une division 
harmonique, si l'on permute ces points, le birapport 
peut prendre trois valeurs, qui sont $-1$, $2$ et $1/2$. 
Il y a donc $4!/3=8$ permutations de $(p_1,p_2,p_3,p_4)$
qui conservent le birapport.
D'autre part, le groupe $\D_8$ op\`ere sur $\{p_1,p_2,p_3,p_4\}$
et le noyau de l'op\'eration est $\{\pm {\rm I}\}$. 
On en d\'eduit que les $8$ permutations ci-dessus sont r\'ealisables 
avec des \'el\'ements de $\D_8$. On se ram\`ene ainsi au cas o\`u $\phi$ 
conserve chaque feuilletage de $\cal{T}_0$ : $\phi$ est une 
dilatatation d'apr\`es le lemme pr\'ec\'edent.
La deuxi\`eme partie de l'\'enonc\'e r\'esulte aussi du lemme pr\'ec\'edent.
\epf

\sk
Rappelons que si $u$ est une fonction d\'efinie au voisinage de $(x_0,y_0)$
qui v\'erifie la condition (de transversalit\'e) 
\beq
\label{transverse}
u_x(x_0,y_0)u_y(x_0,y_0)(u_x(x_0,y_0)^2-u_y(x_0,y_0)^2) \neq 0,
\eeq
($u_x$ et $u_y$ d\'esignent les d\'eriv\'ees partielles de u),
on a d\'efini :
$$
\cal{T}[u]: = \cal{T}(\, x,\, y, \, x+y,\, x-y, u(x,y) \,).
$$
\ble
\label{LC}
Si $u(x,y)$ et $u'(x,y)$ sont deux germes de fonctions qui v\'erifient
(\ref{transverse}) et si $\phi$ un diff\'eomorphisme local
qui envoie $\cal{T}[u]$ sur $\cal{T}[u']$, $\phi$ est 
une transformation projective. Si de plus $\cal{F}(u)$ n'est pas un faisceau
de droites parall\`eles, $\phi$ est un \'el\'ement de ${\rm Sym}\,(\cal{T}_0)$.
\ele
\bpf
Si $\phi$ est comme dans l'\'enonc\'e, au moins trois des faisceaux
de droites de $\cal{T}_0$ sont envoy\'es sur des faisceaux de droites 
de $\cal{T}_0$. On en d\'eduit qu'il existe une transformation 
affine $\psi$ telle $\psi\circ \phi$ conserve ces trois faisceaux.
D'apr\`es le Lemme~\ref{LA}, $\phi$ est une transformation affine.
Si $\cal{F}(u)$ n'est pas un faisceau de droites
parall\`eles, $\phi$ conserve n\'ecessairement le
sous-tissu $\cal{T}_0$ ; le Lemme \ref{LB} s'applique.
\epf
\ble
\label{LD}
Si $u(x,y)$ est une fonction qui v\'erifie (\ref{transverse}), 
le tissu $\cal{T}[u]$ est alg\'ebrisable si et seulement s'il est alg\'ebrique, 
si et seulement si $\cal{F}(u)$ est un germe de faisceau de droites.
\ele
\bpf
On sait (voir l'introduction) qu'un tissu form\'e de faisceaux de droites est alg\'ebrique. 
On suppose que le germe de tissu $\cal{T}[u]$ est alg\'ebrisable. Soit 
$\phi$ un diff\'eomorphisme local tel que le tissu $\phi(\cal{T}[u])$
soit un germe de tissu alg\'ebrique. Il est associ\'e 
\`a une courbe r\'eduite $C_5$ de $\cP$, de degr\'e $5$.

En particulier, $\phi(\cal{T}[u])$ est un tissu en droites, donc 
son sous-tissu $\phi(\cal{T}_0)$ est un tissu en droites. De plus 
$\phi(\cal{T}_0)$ est de rang maximal. C'est donc,
d'apr\`es le Th\'eor\`eme~\ref{abel-inv}, un tissu alg\'ebrique, 
associ\'e \`a une courbe de degr\'e $4$ $C_4\subset C_5$. 
Il est alors clair que $C_5$ est la r\'eunion de $C_4$ et 
d'une droite, donc que $\cal{F}(u)$ est un faisceau de droites.
En appliquant le m\^eme argument \`a tous les $3$-tissus
extraits de $\cal{T}_0$ (ils sont de rang maximal), on montre 
que la courbe $C_4$ est une r\'eunion de droites. Finalement,
$\phi(\cal{T}[u])$ est form\'e de cinq faisceaux de droites,
donc $\phi$ est une transformation projective 
d'apr\`es le Lemme \ref{LA} et $\cal{T}[u]$ est un tissu 
alg\'ebrique.
\epf

\sct{Fonctions th\^eta et tissus exceptionnels}

Cette section contient nos principaux r\'esultats. Pour la plus grand part, 
elle peut \^etre lue ind\'ependamment des autres sections. Nous 
pr\'esentons une famille \`a un param\`etre de $5$-tissus exceptionnels,
ainsi que cinq $5$-tissus exceptionnels \og limites \fg.
Tous ces tissus sont construits par l'adjonction d'un feuilletage 
$\cal{F}(u)$ au tissu $\cal{T}_0= \cal{T}( \, x, \, y, \, x+y, \, x-y \,)$, {\em i.e.}
sont de la forme (\ref{Tu}) :
$$
\cal{T}[u]:= \cal{T}( \, x, \, y, \, x+y, \, x-y, \, u(x,y) \,),
$$
o\`u $u$ v\'erifie la condition de transversalit\'e (\ref{transverse}). 
De plus, ils ont tous la propri\'et\'e que le sous-tissu 
$$
\cal{T}(\, x, \, y, \, u(x,y))
$$
est de rang maximal. 

{\em Dans la section \S 4, nous montrerons que tout 
tissu exceptionnel de la forme $\cal{T}[u]$, tel que 
$\cal{T}(\, x, \, y, \, u(x,y) \,)$ est de rang maximal, est \'equivalent 
\`a l'un des tissus pr\'esent\'es dans cette section.}

D'autre part, mais on verra dans la section \S 4 que c'est
une cons\'equence des conditions ci-dessus,
ces tissus ont la propri\'et\'e que le sous-tissu 
$$
\cal{T}(\, x+y, \, x-y, \, u(x,y))
$$
est aussi de rang maximal. 

\sk
Un tissu de la forme $\cal{T}[u]$ poss\`edent les relations ab\'eliennes du tissu 
$\cal{T}_0$, dont voici une base :
\begin{align*}
x     +  y     - (x+y)                 &= 0,     \\
x     -  y     - (x-y)                 &= 0,     \\
2x^2  +  2y^2  - (x+y)^2  -   (x-y)^2  &=0.
\end{align*}
Le tissu $\cal{T}[u]$ est de rang maximal $\rho(5)=6$ si et seulement s'il
poss\`ede trois relations ab\'eliennes suppl\'ementaires, ind\'ependantes 
modulo les relations pr\'ec\'edentes et les constantes, autrement dit trois 
relations de la forme :
$$
f_k(u(x,y)) = g_k(x) + h_k(y) + j_k(x+y) + l_k(x-y), \qquad k=1,2,3,
$$
avec $f'_1,f'_2,f'_3$ lin\'eairement ind\'ependantes au point $u(x_0,y_0)$. 
D'autre part, les fonctions $u(x,y)$ seront de la forme :
\beq
\label{formeadd}
u(x,y) = v(x) + w(y).
\eeq
C'est une relation ab\'elienne du sous-tissu 
$\cal{T}(x,y,u(x,y))$. Pour prouver que le rang 
est maximal, il restera encore, dans chaque cas,
\`a exhiber deux relations ab\'eliennes. La condition d'ind\'ependance sera 
facile \`a v\'erifier dans tous les cas. 
\bre
Souvent, pour faire l'\'economie d'un logarithme, on \'ecrira 
les relations ab\'eliennes  sous forme multiplicative. En particulier, 
on pourra choisir une fonction de d\'efinition de la forme :
\beq
\label{formemul}
u(x,y) = v(x) w(y).
\eeq
\ere

\bk
On note $\cal{H}:=\{\tau \in \C, \; {\rm Im}\, \tau >0\}$ le demi-plan de Poincar\'e.
Le groupe des homographies $\tau \mapsto (a\tau + b)/(c\tau +d)$
\`a coefficients entiers et de d\'eterminant $1$, ou groupe modulaire, op\`ere 
sur $\cal{H}$. 

\sk
Si $\tau \in \cal{H}$, on note $q:= \e{i\pi \tau}$ et on associe \`a $\tau$ les quatre fonctions th\^eta :
\begin{align*}
\theta_1(x,\tau) & = -i\sum_{n=-\infty}^{+\infty} (-1)^n q^{(n+1/2)^2} \e{i(2n+1)x}, \\   
\theta_2(x,\tau) & =   \sum_{n=-\infty}^{+\infty}        q^{(n+1/2)^2} \e{i(2n+1)x}, \\
\theta_3(x,\tau) & =   \sum_{n=-\infty}^{+\infty}        q^{n^2}      \e{i2nx},      \\
\theta_4(x,\tau) & =   \sum_{n=-\infty}^{+\infty} (-1)^n q^{n^2}      \e{i2nx}.
\end{align*}
Pour all\'eger les formules, on notera simplement :
$$
\theta_i(x):= \theta_i(x,\tau), \qquad i=1,2,3,4,
$$
quand la valeur du param\`etre est $\tau$ ; sinon, on \'ecrira le param\`etre.

Les fonctions th\^eta sont des fonctions enti\`eres de $x\in \C$. La fonction $\theta_1$ est impaire 
et les fonctions $\theta_2$, $\theta_3$ et $\theta_4$ sont paires. 
La formule fondamentale suivante est classique
et se d\'emontre facilement \`a partir des d\'efinitions ci-dessus (voir par exemple
Lawden \cite{La}, chapitre 1) :
$$
\theta_3(x) \, \theta_4(y) 
=
\theta_4(x+y,2\tau) \, \theta_4(x-y,2\tau) - \theta_1(x+y,2\tau) \, \theta_1(x-y,2\tau).
$$
Apr\`es changement de variables :
\beq
\label{form}
\theta_3((x+y)/2,\tau/2) \,   \theta_4((x-y)/2,\tau/2) 
= 
\theta_4(x) \,  \theta_4(y) - \theta_1(x) \,  \theta_1(y).
\eeq
On a le r\'esultat suivant :
\bt
\label{famille}
Pour tout $\tau\in \cal{H}$, le tissu $\cal{T}[u]$ associ\'e \`a la fonction 
$$
u(x,y) = \fr{\theta_1(x)}{\theta_4(x)} \, \fr{\theta_1(y)}{\theta_4(y)},
$$
en un point o\`u (\ref{transverse}) est v\'erifi\'e,
est exceptionnel.
Les tissus associ\'es \`a $\tau,\tau' \in \cal{H}$
sont \'equivalents si et seulement si $\tau$ et $\tau'$ sont congrus 
modulo le groupe $G$ d'automorphismes de $\cal{H}$
engendr\'e par les transformations $\tau \mapsto \tau+2$, $\tau \mapsto \tau/(\tau + 1)$
et $\tau\mapsto -2/\tau$.
\et
Notons que $\tau \mapsto -2/\tau$ n'est pas une transformation modulaire.
La deuxi\`eme partie de l'\'enonc\'e, qu'on d\'emontrera \`a la fin de cette section,
a le sens suivant : si $\tau$ et $\tau'\in \cal{H}$ sont congrus modulo 
le groupe $G$, les tissus associ\'es sont \'equivalents par un \'el\'ement du 
groupe ${\rm Sym}\,(\cal{T}_0)$ ; de plus (c'est une cons\'equence 
de la propri\'et\'e pr\'ec\'edente et du Lemme \ref{LC}), si des germes des
tissus associ\'es \`a $\tau, \tau'\in \cal{H}$, en des points 
o\`u la condition (\ref{transverse}) est v\'erifi\'ee, sont \'equivalents,
alors $\tau$ et $\tau'$ sont congrus modulo $G$.
\bpf
On obtient une relation ab\'elienne (sous forme multiplicative)
en divisant les deux membres de (\ref{form}) 
par $\theta_4(x) \,  \theta_4(y)$ :
$$
1 - u(x,y) =  
\fr{\theta_3((x+y)/2,\tau/2) \, \theta_4((x-y)/2,\tau/2)}
   {\theta_4(x)\,\theta_4(y)}
$$
et une autre en rempla\c{c}ant $(x,y)$ par $(x,-y)$ 
et en tenant compte des parit\'es des fonctions th\^eta :
$$
1 + u(x,y) =  
\fr{\theta_3((x-y)/2,\tau/2) \, \theta_4((x+y)/2,\tau/2)}
   {\theta_4(x)\,\theta_4(y)}.
$$
On v\'erifie que $\cal{F}(u)$ n'est pas un faisceau de droites !
Compte tenu du Lemme \ref{LD}, le tissu $\cal{T}[u]$ est exceptionnel.
\epf

\bk
En divisant membre \`a membre les deux relations pr\'ec\'edentes, on obtient une relation \`a $3$ facteurs :
\beq
\label{e4}
\fr{1 - u(x,y)}{1 + u(x,y)}
=  
\fr{\theta_3((x+y)/2,\tau/2)/\theta_4((x+y)/2,\tau/2)}
   {\theta_3((x-y)/2,\tau/2)/\theta_4((x-y)/2,\tau/2)}.
\eeq
Plus bas, nous rappellerons la d\'efinition des fonctions elliptiques 
de Jacobi. La relation (\ref{e4}) s'\'ecrit simplement en termes de ces 
fonctions. 

\sk
Pour faire la liaison entre cette famille et les cinq
tissus exceptionnels que nous allons pr\'esenter maintenant, 
notons seulement que la famille de tissus du Th\'eor\`eme~\ref{famille}
est \'equivalente \`a la famille $\cal{T}[u_k]$ d\'efinie par les fonctions 
$$
k\in \C, \; k^2\neq 0,1, \qquad  u_k(x,y) = {\rm sn}_k x \, {\rm sn}_k y,
$$
et qu'on a les limites classiques suivantes :
$$
\tanh x = \lim_{k\rightarrow 1} {\rm sn}_k x, 
\qquad
\sin x =  \lim_{k\rightarrow 0} {\rm sn}_k x.
$$

\bt
\label{autres}
On obtient cinq $5$-tissus exceptionnels $\cal{T}[u]$ deux \`a deux non \'equivalents 
pour les choix suivants de la fonction $u(x,y)$ :
\be
\item[] \hspace{0.5cm} {\rm (A)} \hspace{2cm}  $u(x,y) = \tanh x \tanh y$,
\item[] \hspace{0.5cm} {\rm (B)} \hspace{2cm}  $u(x,y) = \sin x \sin y$, 
\item[] \hspace{0.5cm} {\rm (C)} \hspace{2cm}  $u(x,y) = \e{x} + \e{y}$, 
\item[] \hspace{0.5cm} {\rm (D)} \hspace{2cm}  $u(x,y) = x^2 - y^2$, 
\item[] \hspace{0.5cm} {\rm (E)} \hspace{2cm}  $u(x,y) = x^2 + y^2$. 
\ee
\et
\bre
Les sym\'etries $(x,y)\mapsto (\pm x,\pm y)$ et $(x,y) \mapsto (\pm y, \pm x)$
du tissu $\cal{T}_0$ transforment les tissus (A)--(D) en des tissus 
qui leur sont \'equivalents par dilatation. En d\'emontrant le th\'eor\`eme, 
on verra que la sym\'etrie $\sigma(x,y)=((x+y)/\sqrt{2},(x-y)/\sqrt{2})$ 
transforme les tissus 
(A)--(D) en des tissus de la forme $\cal{T}[u_\sigma]$ qui, modulo une dilatation,
peuvent \^etre d\'efinis par les fonctions suivantes :
\be
\item[] \hspace{0.5cm} {\rm (A)} \hspace{2cm}  $u_\sigma(x,y) = \cosh x/\cosh y$,
\item[] \hspace{0.5cm} {\rm (B)} \hspace{2cm}  $u_\sigma(x,y) = \cos x + \cos y$, 
\item[] \hspace{0.5cm} {\rm (C)} \hspace{2cm}  $u_\sigma(x,y) = \e{x}\cosh y$, 
\item[] \hspace{0.5cm} {\rm (D)} \hspace{2cm}  $u_\sigma(x,y) = xy$, 
\item[] \hspace{0.5cm} {\rm (E)} \hspace{2cm}  $u_\sigma(x,y) = x^2+y^2$. 
\ee
\ere
\bre
La remarque qui pr\'ec\`ede l'\'enonc\'e montre que les tissus (A) et (B) sont 
\og des cas limites \fg\, des tissus du Th\'eor\`eme \ref{famille}.
Les limites ci-dessous montrent qu'on peut consid\'erer les tissus (C), (D) et (E) 
comme des limites de suites de tissus tous \'equivalents au tissu (B). 
\begin{align*}
\e{x} + \e{y} & = \lim_{k\rightarrow +\infty} 2\e{-k}(\cosh (x+k) + \cosh (y+k)), \\
x^2 \pm  y^2  & = \; \lim_{\epsilon \rightarrow 0} \;\,
2\epsilon^{-2}((\cosh (\epsilon x)-1) \pm (\cosh (\epsilon y)-1)).
\end{align*}
\ere
\bpf
Elle consiste \`a exhiber, dans chaque cas, 
trois relations ab\'eliennes qui, avec les relations
de $\cal{T}_0$, prouvent que le tissu $\cal{T}[u]$ est de rang maximal.
Selon les cas, il est plus commode d'\'ecrire ces relations  
sous forme additive ou sous forme multiplicative.
Ceci fait, le fait que $\cal{T}[u]$ est exceptionnel est une cons\'equence 
imm\'ediate du Lemme~\ref{LD}. 

Pour construire une fonction de d\'efinition $u_\sigma$ 
de la forme (\ref{formeadd}) ou (\ref{formemul})
du tissu $\sigma(\cal{T}[u])$, en fait modulo une dilatation, 
on \'ecrit une relation ab\'elienne 
du sous-tissu $\cal{T}(x+y,x-y,u(x,y))$ et on la transforme par $\sigma$.

La non \'equivalence des cinq tissus 
se v\'erifie, compte tenu du Lemme \ref{LC}, en comparant les dix 
fonctions $u$ ou $u_\sigma$ obtenues.

Dans chaque cas, on indiquera le cas \'ech\'eant un changement de variables qui 
transforme le tissu $\cal{T}[u]$ en un tissu dont les feuilles sont des
(germes de) courbes alg\'ebriques. 

\bk
{\bf (A). --- } On a les relations ab\'eliennes : 
\begin{align*}
u(x,y)         & = \tanh x \tanh y ,                 \\
1 + u(x,y)     & = \fr{\cosh (x+y)}{\cosh x\cosh y}, \\
1 - u(x,y)     & = \fr{\cosh (x-y)}{\cosh x\cosh y}.
\end{align*}
En divisant membre \`a membre les deux derni\`eres identit\'es, on obtient :
$$
\fr{1 + u(x,y)}{1 - u(x,y)} = \fr{\cosh (x + y) }{\cosh (x - y)},
$$
dont on d\'eduit :
$$
u_\sigma (x,y) = \cosh x/\cosh y.
$$
En posant $\e{2x}=\xi$, $\e{2y}=\eta$, on obtient que le tissu $\cal{T}[\, \tanh x \, \tanh y \, ]$
est \'equivalent au tissu 
$$
\cal{T}(\, \xi, \, \eta, \, \xi \eta, \, \xi/\eta, \, \fr{(\xi - 1)(\eta - 1)}{(\xi + 1)(\eta + 1)}\,),
$$
form\'e de trois faisceaux de droites et de deux faisceaux de coniques.

\bk
{\bf (B). --- } On a les relations ab\'eliennes :
\begin{align*}
 u(x,y)      & = \sin x \sin y,                               \\
2u(x,y)      & = \cos(x-y)   - \cos(x+y),                     \\
4u(x,y)^2    & = \cos^2(x-y) + \cos^2(x+y) - \cos 2x - \cos 2y.
\end{align*}
De la deuxi\`eme relation, on d\'eduit :
$$
u_\sigma(x,y) = \cos x  + \cos y.
$$
En posant $\e{ix}=\xi$, $\e{iy}=\eta$, on obtient que le tissu $\cal{T}[\, \sin x \, \sin y \,]$
est \'equivalent au tissu 
$$
\cal{T}(\, \xi, \, \eta, \, \xi \eta, \, \xi/\eta, \, \fr{(\xi^2 - 1)(\eta^2 - 1)}{\xi \eta}\, ),
$$
form\'e de trois faisceaux de droites, un faiseau de coniques et un faisceau de quartiques.

\bk
{\bf (C). --- } On a les relations ab\'eliennes :
\begin{align*}
u(x,y)     & = \e{x} + \e{y},                     \\
u(x,y)^2   & = \e{2x} + \e{2y} + 2\e{(x+y)},        \\
u(x,y)     & = 2\e{(x+y)/2} \cosh \fr{x-y}{2}.
\end{align*}
De la derni\`ere relation, on d\'eduit :
$$
u_\sigma(x,y)  = \e{x} \cosh y.
$$
En posant $\e{x}=\xi$, $\e{y}=\eta$, on obtient que le tissu $\cal{T}[\, \e{x}+\e{y}\, ]$
est \'equivalent au tissu 
$$
\cal{T}(\, \xi, \, \eta, \, \xi \eta, \, \xi/\eta, \, \xi + \eta\, ),
$$
form\'e de quatre faisceaux de droites et d'un faisceau de coniques.

\bk
{\bf (D). --- } On a les relations ab\'eliennes :
\begin{align*}
u(x,y)         & = x^2 - y^2,                            \\
6u(x,y)^2      & = 8x^4 + 8y^4 - (x+y)^4 - (x-y)^4,      \\
u(x,y)         & = (x+y)(x-y).
\end{align*}
De la derni\`ere relation, on d\'eduit :
$$
u_\sigma(x,y) = xy.
$$

\bk
{\bf (E). --- } On a les relations ab\'eliennes :
\begin{align*}
u(x,y)        & = x^2  + y^2,                           \\
6u(x,y)^2     & = 4x^4  + 4y^4  + (x+y)^4  + (x-y)^4,     \\
10u(x,y)^3    & = 8x^6 + 8y^6 + (x+y)^6  + (x-y)^6.
\end{align*}
La relation \`a trois termes $2u(x,y)^2 = (x+y)^2 + (x-y)^2$
est (aussi) une relation du $4$-tissu $\cal{T}_0$. 
Le tissu $\cal{T}[x^2+y^2]$ est invariant par $\sigma$.
\epf

On va maintenant traduire le Th\'eor\`eme \ref{famille}
en termes des fonctions elliptiques de Jacobi, en d\'emontrer la
deuxi\`eme partie et donner des mod\`eles de ses tissus 
qui ont la propri\'et\'e d'\^etre \og \`a feuilles alg\'ebriques\fg.

\sk
On note :
$$
\theta_i:= \theta_i(0,\tau), \qquad i=1,2,3,4. 
$$
\`A $\tau \in \cal{H}$ on associe les nombres :
$$
k = \theta_2^2/\theta_3^2, \qquad  k' = \theta_4^2/\theta_3^2.
$$
Les fonctions elliptiques {\em de param\`etre $\tau$} ou {\em de module $k$} sont d\'efinies par :
$$
{\rm sn}_k x = \fr{\theta_3}{\theta_2} \fr{\theta_1(x/\theta_3^2)}{\theta_4(x/\theta_3^2)}, 
\;\;
{\rm cn}_k x = \fr{\theta_4}{\theta_2} \fr{\theta_2(x/\theta_3^2)}{\theta_4(x/\theta_3^2)},
\;\;
{\rm dn}_k x = \fr{\theta_4}{\theta_3} \fr{\theta_3(x/\theta_3^2)}{\theta_4(x/\theta_3^2)}.
$$
Le nombre $k'$ est appel\'e {\em le module conjugu\'e de $k$}. On a $k^2 + k'^2=1$.
Ces fonctions ne d\'ependent en fait que du carr\'e du module. La fonction 
\beq
\label{modulaire}
\tau \in \cal{H}, \qquad \tau\mapsto k(\tau)^2
\eeq
est une fonction modulaire, c'est-\`a-dire 
invariante par un sous-groupe d'indice fini du groupe modulaire. Elle prend 
toutes les valeurs complexes sauf $0$ et $1$. Elle est \'etudi\'ee 
par exemple dans Chandrasekharan \cite{Ch}, livre auquel 
on renvoie pour plus de d\'etails, en particulier sur les relations 
entre le param\`etre $\tau$ et le module $k$, qu'on ne discutera
pas. 

Avec ces notations, la famille de fonctions $u(x,y)$ du Th\'eor\`eme \ref{famille}
s'\'ecrit $u(x,y) =  k\, {\rm sn}_k(\theta_3^2x)\, {\rm sn}_k(\theta_3^2y)$.
Modulo les dilatations de $\C^2$, on peut remplacer cette famille 
par la suivante :
\beq
\label{form-elli}
k\in \C, \; k^2\neq 0,1, \qquad u_k(x,y) = {\rm sn}_k x \,\, {\rm sn}_k y.
\eeq
Compte tenu du Lemme \ref{LD} et de propri\'et\'es de la fonction (\ref{modulaire})
qu'on trouvera dans \cite{Ch}, la deuxi\`eme partie du Th\'eor\`eme \ref{famille} 
est une cons\'equence de la proposition suivante :
\bpr
Soit $k,l\in \C$ avec $k^2,l^2\neq 0,1$. Les tissus (avec singularit\'es)
$\cal{T}[\,{\rm sn}_k x \,\, {\rm sn}_k y\,]$ et $\cal{T}[\,{\rm sn}_l x \,\, {\rm sn}_l y\,]$ 
sont \'equivalents sous l'action du groupe ${\rm Sym}\,(\cal{T}_0)$ si et seulement si
\beq
\label{ketl}
l \in \{ \, \pm k, \, \pm \fr{1}{k}, \, \pm \fr{1-k}{1+k}, \, \pm \fr{1+k}{1-k} \, \}.
\eeq
\epr
Compte tenu du Lemme \ref{LC}, cette proposition r\`egle aussi la question 
de l'\'equivalence dans les germes.
Avant de d\'emontrer cette proposition, nous \'ecrivons la relation (\ref{e4}) 
sous une autre forme :
\ble
\label{lemme-e4bis}
On a l'identit\'e :
\beq
\label{e4bis}
\fr{1 + k\,{\rm sn}_k x \, {\rm sn}_k y}{1 - k\,{\rm sn}_k x \, {\rm sn}_k y}
 = 
\fr{{\rm dn}_k (x+y) - k\,{\rm cn}_k (x+y)}{{\rm dn}_k (x-y) - k\,{\rm cn}_k (x-y)}.
\eeq
\ele
\bpf
On pourrait partir de (\ref{e4}) mais on utilisera plut\^ot 
le formulaire des fonctions elliptiques. On part des formules d'additions
(voir \cite{La} (2.4.1)--(2.4.3) ; on n'\'ecrit pas l'indice $k$) :
\begin{align*}
\cn (x\pm y)  & = \fr{\cn x \, \cn y     \mp     \, \sn x \, \sn y \, \dn x \, \dn y}{1 - k^2 \, \sn^2 x \, \sn^2 y},
\\
\dn (x\pm y)  & = \fr{\dn x \, \dn y     \mp k^2 \, \sn x \, \sn y \, \cn x \, \cn y}{1 - k^2 \, \sn^2 x \, \sn^2 y}.
\end{align*}
En substituant les seconds membres 
\`a $\cn (x\pm y)$ et $\dn (x\pm y)$ dans le second membre de (\ref{e4bis}), on obtient le r\'esultat.
\epf
\bre
En rempla\c{c}ant $(x,y)$ par $\sigma(x,y)$ dans (\ref{e4bis}), on voit 
que le tissu $\sigma(\cal{T}[\,{\rm sn}_k x \,\, {\rm sn}_k y\,])$ est 
un dilat\'e du tissu $\cal{T}[v]$, 
o\`u $v(x,y)=({\rm dn}_k x - k{\rm cn}_k x)/({\rm dn}_k y - k{\rm cn}_k y)$.
\ere

\sk
Pour \'eviter d'\'ecrire des formules pr\'ecises de la th\'eorie des fonctions elliptiques, 
on convient de noter :
\beq
\label{dilat}
u'(x,y) \sim u(x,y) 
\eeq
s'il existe $c,k\in \C^*$ et $x',y',u_0\in \C$ tels qu'on ait l'identit\'e
$$
(x,y)\in \C^2, \qquad u'(x,y) = cu(kx+x',ky+y') + u_0.
$$
On utilisera une convention analogue pour les fonctions $u(x)$ d'une variable.
\bpf[D\'emonstration de la Proposition]
Pour un module $k$ fix\'e, on a les formules \og de p\'eriodicit\'e\fg\, suivantes,
voir (\cite{La} (2.1.21), (2.2.17)) :
\beq
\label{periode}
{\rm sn}_k x  \sim 1/{\rm sn}_k x,  \qquad {\rm sn}_k x  \sim {\rm dn}_k x/{\rm cn}_k x.
\eeq
Par exception, on aura  besoin de la version plus pr\'ecise suivante de la premi\`ere formule :
il existe $T\in \C$, qui ne d\'epend que de $k$, tel que :
\beq
\label{inverse-sn}
x\in \C, \qquad k \, {\rm sn}_k (x+T) =  1/{\rm sn}_k x.
\eeq
D'autre part, on a les formules suivantes, qui relient des fonctions elliptiques 
de modules diff\'erents (voir \cite{La} (3.9.4), (3.9.24), (2.6.12)) :
$$
{\rm sn}_{1/k} x              \sim {\rm sn}_k x,        \qquad 
{\rm dn}_k x - k{\rm cn}_k x  \sim 1/{\rm dn}_{l'} x,   \qquad
{\rm dn}_{k'} x               \sim  {\rm dn}_k x/ {\rm cn}_k x,
$$
o\`u $l=(1-k)/(1+k)$ et $k'$ et $l'$ sont les modules conjugu\'es de $k$ et $l$.
On a donc :
$$
\fr{{\rm dn}_k x - k{\rm cn}_k x}{{\rm dn}_k y - k{\rm cn}_k y}
\sim
\fr{{\rm dn}_{l'} y}{{\rm dn}_{l'} x}
\sim
\fr{{\rm dn}_l y/{\rm cn}_l y}{{\rm dn}_l x/{\rm cn}_l x}
\sim
\fr{{\rm sn}_l y}{{\rm sn}_l x}
\sim
{\rm sn}_l y \, {\rm sn}_l x.
$$
Il r\'esulte du Lemme \ref{lemme-e4bis}, de la remarque qui le suit 
et des \'equivalences ci-dessus que 
que le tissu $\cal{T}[\,{\rm sn}_l x \,\, {\rm sn}_l y\,]$ est \'equivalent 
par dilatation \`a l'image du tissu 
$\cal{T}[\,{\rm sn}_k x \,\, {\rm sn}_k y\,]$
par la sym\'etrie $\sigma$.
Avec (\ref{periode}), ceci montre que, dans l'\'enonc\'e de la proposition, la condition (\ref{ketl})
est suffisante. 

R\'eciproquement, on suppose que 
$\cal{T}[\,{\rm sn}_k x \,\, {\rm sn}_k y\,]$ et 
$\cal{T}[\,{\rm sn}_l x \,\, {\rm sn}_l y\,]$ 
sont dans la m\^eme orbite sous l'action du groupe ${\rm Sym}\,(\cal{T}_0)$.
Quitte \`a remplacer $l$ par $(1-l)/(1+l)$, on peut 
supposer qu'ils sont dans la m\^eme orbite sous l'action 
du groupe engendr\'e par les dilatations de $\C^2$
et les sym\'etries $(x,y)\mapsto (\pm x, \pm y)$ et 
$(x,y)\mapsto (\pm y, \pm x)$ (\og on a enlev\'e $\sigma$ \fg).
Comme ces sym\'etries laissent ces tissus invariants, 
on peut supposer qu'ils sont \'equivalents 
par dilatation. On a alors n\'ecessairement 
${\rm sn}_k x \,\, {\rm sn}_k y   \sim {\rm sn}_l x \,\, {\rm sn}_l y$, donc :
$$
x\in \C, \qquad {\rm sn}_k x =  c\,{\rm sn}_l (dx+x') + e,
$$
avec $c,d\in \C^*$ et $x',e\in \C$. D'autre part, et c'est une propri\'et\'e classique,
la fonction ${\rm sn}_k x$ est solution de l'\'equation diff\'erentielle :
$$
(z')^2 = (1-z^2)(1-k^2z^2).
$$
On en d\'eduit que ${\rm sn}_l x$ v\'erifie l'\'equation diff\'erentielle :
$$
(z')^2 = (1-(cz+e)^2)(1-k^2(cz+e)^2)/(c^2d^2).
$$
Le second membre doit \^etre identique au second membre de 
l'\'equation diff\'erentielle $(z')^2 = (1-z^2)(1-l^2z^2)$. On 
obtient $e=0$, puis $l^2=k^2$ ou $l^2=1/k^2$. Ceci termine 
la d\'emonstration de la proposition.
\epf

\sk
On obtient un \og mod\`ele \`a feuilles alg\'ebriques \fg\, 
du tissu $\cal{T}[{\rm sn}_k x \, {\rm sn}_k y]$ en posant :
$$
{\rm sn}_k^2 x = \xi, \qquad {\rm sn}_k^2 y = \eta. 
$$
On obtient le tissu \'equivalent :
$$
\cal{T}(\, \xi, \, \eta, \, u_+(\xi,\eta), \, u_-(\xi,\eta), \, \xi \eta \, ),
$$
o\`u $u_\pm (\xi,\eta) = {\rm sn}_k^2 (x\pm y)$ peut \^etre calcul\'e gr\^ace aux formules
d'additions 
$$
{\rm sn}_k (x\pm y) = 
\fr{{\rm sn}_k x \, {\rm cn}_k y \, {\rm dn}_k y \, \pm  \,{\rm sn}_k y \, {\rm cn}_k x \, {\rm dn}_k x}
{1 - k^2 \, {\rm sn}_k^2 x \, {\rm sn}_k^2 y}
$$
et aux identit\'es :
$$
{\rm sn}_k^2 x + {\rm cn}_k^2 x = 1, \qquad  {\rm dn}_k^2 x + k^2\, {\rm sn}_k^2 x = 1.
$$
On \'ecrit :
$$
{\rm sn}_k(x\pm y) = \sqrt{A(\xi,\eta)} \pm \sqrt{A(\eta,\xi)},
$$
avec 
$$
A(\xi,\eta) = \fr{{\rm sn}^2_k x \, {\rm cn}^2_k y \, {\rm dn}^2_k y}{(1 - k^2\, {\rm sn}_k^2 x \, {\rm sn}_k^2 y)^2}
     = \fr{\xi (1-\eta)(1-k^2\eta)}{(1 - k^2 \xi\eta)^2}.
$$
On a donc :
\begin{align*}
u_+(\xi,\eta) + u_-(\xi,\eta) & = 2(A(\xi,\eta)+A(\eta,\xi)), \\
u_+(\xi,\eta)u_-(\xi,\eta)    & = (A(\xi,\eta)-A(\eta,\xi))^2.
\end{align*}
On v\'erifie que 
$$
A(\xi,\eta)-A(\eta,\xi) = \fr{\xi - \eta}{1-k^2\xi\eta}.
$$
Au point $(\xi,\eta)$ de $\C^2$, $u_+(\xi,\eta)$ et $u_-(\xi,\eta)$ sont les deux racines de l'\'equation 
en $t$ :
\beq
\label{quartique}
(1-k^2\xi\eta)^2t^2 - 2(\xi(1-\eta)(1-k^2\eta)+\eta(1-\xi)(1-k^2\xi)t + (\xi-\eta)^2 = 0.
\eeq
On en d\'eduit que les feuilles des feuilletages $\cal{F}(u_+)$ et $\cal{F}(u_-)$
sont des germes des quartiques obtenues en fixant $t\in \C$
dans l'\'equation (\ref{quartique}).

\sct{Un syst\`eme diff\'erentiel associ\'e aux tissus de la section \S 3}

Dans cette section, nous d\'emontrons le r\'esultat suivant :
\bt
\label{tous}
Soit $u$ un germe de fonction qui v\'erifie (\ref{transverse}). 
Si le tissu $\cal{T}[u]$ est exceptionnel 
et si le tissu $\{\cal{F}(x), \cal{F}(y), \cal{F}(u)\}$ est de rang maximal,
le tissu $\cal{T}[u]$ est l'image, par un \'el\'ement du groupe 
${\rm Sym}\,(\cal{T}_0)$, d'un germe d'un des tissus consid\'er\'es dans 
les  Th\'eor\`emes \ref{famille} et \ref{autres}.
\et
La d\'emonstration repose sur la simplicit\'e du syst\`eme diff\'erentiel 
que doit v\'erifier une fonction $u(x,y)$, choisie de la forme $v(x)+w(y)$,
pour que le tissu $\cal{T}[u]$ soit de rang maximal.

Plus g\'en\'eralement, on peut obtenir une condition diff\'erentielle 
assez simple (\`a \'ecrire, pas toujours \`a r\'esoudre)
sur la fonction $u(x,y)$ pour qu'un tissu de la forme 
\beq
\label{gene}
\cal{T}( \, l_1, \, \ldots, \, l_p, \, u \,),
\eeq
o\`u $l_1,\ldots,l_p$ sont des formes lin\'eaires deux \`a deux lin\'eairement 
ind\'e\-pen\-dantes, soit de rang maximal. Les feuilles des feuil\-letages 
qui constituent le tissu (\ref{gene}) sont aussi les 
les courbes int\'egrales des champs de vecteurs 
$$
X_k = \fr{\pl l_k}{\pl y}\fr{\pl}{\pl x} - \fr{\pl l_k}{\pl x}\fr{\pl }{\pl y},
\;\; k=1,\ldots,p \, ; \qquad
X_u = \fr{\pl u}{\pl y}\fr{\pl}{\pl x} - \fr{\pl u}{\pl x}\fr{\pl }{\pl y}.
$$
{\em Les champs $X_1,\ldots,X_p$ sont constants}. Le $p$-tissu
$\cal{T}(l_1,\ldots,l_p)$, form\'e de 
$p$ faisceaux de droites parall\`eles, est de rang maximal 
$\rho(p)=(p-1)(p-2)/2$. En raisonnant comme au d\'ebut de la section \S 3,
on voit  que le $(p+1)$-tissu (\ref{gene}) est de rang maximal
$\rho(p+1) = \rho(p) + p-1$ si et seulement si l'espace 
des fonctions $f'(u)$ telles que $f(u)$
v\'erifie une relation de la forme :
\beq
\label{relgene}
f(u) = \sum_{k=1}^p f_k(l_k)
\eeq
est de dimension $p-1$, si et seulement si 
l'espace des fonctions $f(u)$ de cette forme est de dimension $p$. On a le lemme \'el\'ementaire :
\ble
Avec les notations qu'on vient d'introduire, une fonction 
$h$ est de la forme $h=\sum_{k=1}^p f_k(l_k)$
si et seulement si $X_1\ldots X_p \, h = 0$.
\ele
\bpf
Ce qui, \`a la rigueur, demande une d\'emons\-tration,
est le fait que, si $X_1\ldots X_p \, h = 0$, alors $h$ est de la forme 
annonc\'ee. C'est clair si $p=1$ et, si c'est connu 
\`a l'ordre $p-1$, 
$$
(X_1 \ldots X_{p-1}) \, X_ph = 0
$$
donne $X_ph = \sum_{k=1}^{p-1} f_k(l_k)$. Si, pour $k=1,\ldots p-1$, 
$g_k$ est une primitive de $f_k$, $X_pg_k(l_k)=(X_pl_k)f_k(l_k)$ 
avec $X_pl_k\in \C^*$ (car $l_k$ et $l_p$ sont ind\'ependantes). 
On en d\'eduit l'existence de $g_p$ tel que :
$$
h = \sum_{k=1}^{p-1} g_k(l_k)/(X_pl_k) + g_p(l_p).
$$
\epf
On en d\'eduit le crit\`ere suivant :
\bpr
\label{prop}
Avec les notations pr\'ec\'edentes, soit $a_1,\ldots ,a_p$ 
les coefficients dans la formule 
\beq
\label{equadif}
X_1\ldots X_p \,  f(u) \equiv \sum_{k=1}^p a_k(x,y) f^{(k)}(u),
\eeq
obtenue en calculant formellement $X_1\ldots X_p \, f(u)$.
Pour que le tissu $\cal{T}(\, l_1, \, \ldots, \, l_p, \, u \,)$
soit de rang maximal, il faut et il suffit que $u$ v\'erifie 
le syst\`eme d'\'equations :
\beq
\label{sys}
k=1,\ldots,p-1, \qquad  X_u \, (a_k/a_p) \equiv 0.
\eeq
\epr
On peut calculer facilement les fonctions $a_k$ en termes de $u$ et de 
ses d\'eriv\'ees, mais on n'utilisera que les cas $p=2$ et $p=4$.
Notons seulement que 
$$
a_p \equiv (X_1u)\ldots (X_pu) \neq 0
$$
au point d'\'etude ; c'est la condition de transversalit\'e 
des feuilletages.
\bpf
On note $b_k=a_k/a_p$, $k=1,\ldots,p-1$. Compte tenu de la discussion qui 
pr\'ec\`ede l'\'enonc\'e, 
le tissu est de rang maximal si et seulement si 
l'\'equation 
$$
f^{(p)}(u(x,y)) + \sum_{k=1}^{p-1} b_k(x,y) f^{(k)}(u(x,y)) = 0
$$
a un espace de solutions de dimension $p$. Si l'on prend 
un syst\`eme de coordonn\'ees locales $(u,v)$ dont 
$u$ fait partie, l'\'equation pr\'ec\'edente est 
une \'equation diff\'erentielle ordinaire, avec $u$ comme variable
de d\'erivation et $v$ comme param\`etre. Il est clair 
que le premier membre ne doit pas d\'ependre explicitement de $v$.
D'o\`u la proposition.
\epf

\sk
Revenons au cas particulier des tissus $\cal{T}[u]$, pour lesquels $p=4$ et 
$$
X_1 = \pl/\pl y, \;\; X_2 = \pl/\pl x, \;\; X_3 = \pl/\pl x - \pl/\pl y, \;\; 
X_4 = \pl/\pl x + \pl/\pl y.
$$
On suppose aussi que $u(x,y)$ est de la forme 
$$
u(x,y) = v(x) + w(y)
$$
et v\'erifie la condition (de transversalit\'e) :
\beq
\label{transversebis}
v_xw_y(v_x^2-w_y^2) \neq 0.
\eeq
(Par exemple, $v_x,v_{xx}$... d\'esignent les d\'eriv\'ees successives de $v$). On a :
\begin{align*}
X_1X_2 f(u)      & = (f(u))_{xy} = f''(u)v_xw_y,         \\
X_1X_2X_3X_4f(u) & = (f''(u)v_xw_y)_{xx} - (f''(u)v_xw_y)_{yy},
\end{align*}
avec par exemple :
\begin{align*}
(f''(u)v_xw_y)_{xx} & = w_y(f'''(u)v_x^2 + f''(u)v_{xx})_x  \\
                    & = w_y(f''''(u)v_x^3 +  3f'''(u)v_xv_{xx} + f''(u)v_{xxx}).
\end{align*}
Par sym\'etrie, on obtient :
\begin{multline*}
X_1X_2X_3X_4f(u) = f''''(u)v_xw_y(v_x^2-w_y^2) +  \\
3f'''(u)v_xw_y(v_{xx}-w_{yy}) + f''(u)(w_yv_{xxx}-v_xw_{yyy}).
\end{multline*}
La proposition pr\'ec\'edente s'applique : $u$ d\'efinit un tissu de rang maximal si 
et seulement si les \'equations suivantes sont v\'erifi\'ees :
\begin{align}
\label{eq1}
&  X_u \left(\fr{v_{xx}-w_{yy}}{v_x^2-w_y^2}\right)=0,  \\
\label{eq2}
&  X_u \left(\fr{w_yv_{xxx}-v_xw_{yyy}}{v_xw_y(v_x^2-w_y^2)}\right)=0.
\end{align}
On a le lemme {\em a priori} surprenant suivant :
\ble
Si $u(x,y)=v(x)+w(y)$ v\'erifie (\ref{transversebis}),
l'\'equation (\ref{eq1}) est une condition n\'ecessaire et suffisante pour que le $3$-tissu
$\cal{T}(\, x+y, \, x-y, \, u(x,y)\,)$ soit de rang maximal.
\ele
\bpf
On a :
$$
X_3X_4f(u) = f(u)_{xx} - f(u)_{yy} = f''(u)(v_x^2 - w_y^2) + f'(u)(v_{xx}-w_{yy}).
$$
La Proposition \ref{prop} et la comparaison avec (\ref{eq1})  donnent le r\'esultat.
\epf
On a $X_u=v_x \pl/\pl y - w_y \pl/\pl x$. L'\'equation (\ref{eq1}) s'\'ecrit :
\beq
\label{eq1bis}
v_x \left( \fr{v_{xx}-w_{yy}}{v_x^2-w_y^2}\right )_{\!\!y} = w_y\left( \fr{v_{xx}-w_{yy}}{v_x^2-w_y^2} \right)_{\!\!x}.
\eeq
On est amen\'e \`a traiter \`a part les cas o\`u l'une des fonctions 
$v_{xx}$ ou $w_{yy}$ s'annule identiquement. Si les deux fonctions 
sont nulles, $u$ est une fonction affine et $\cal{F}[u]$ un faisceau de droites 
parall\`eles. Si par exemple $v_{xx}\equiv 0$ et $w_{yy}\not\equiv 0$,
soit $v_x\equiv c \in \C^*$, on obtient :
\beq
\label{special}
\left( \fr{w_{yy}}{c^2-w_y^2}\right )_{\!\!y} = 0,
\eeq
d'o\`u, pour un $k\in \C^*$, $w_{yy} = k(c^2-w_y^2)$ et 
$w_{yyy}=-2kw_yw_{yy}=-2k^2w_y(c^2-w_y^2)$. Pour le moment, on 
se contente de remarquer que l'\'equation (\ref{eq2}), 
qui se r\'eduit dans ce cas \`a :
$$
\left(\fr{w_{yyy}}{w_y(c^2-w_y^2)}\right)_{\!\!y}=0
$$
est impliqu\'ee par l'\'equation (\ref{eq1}).

\sk
{\em On suppose maintenant que les fonctions $v_{xx}$ et $w_{yy}$ ne sont pas 
identiquement nulles}. On d\'eveloppe d'abord (\ref{eq1bis}) sous la forme :
\beq
\label{eq1ter}
( \fr{v_{xxx}}{v_x} + \fr{w_{yyy}}{w_y} )(v_x^2-w_y^2) = 2(v_{xx}^2-w_{yy}^2).
\eeq
Il en r\'esulte que la fonction 
$$
h:= v_x^2\fr{w_{yyy}}{w_y} - w_y^2\fr{v_{xxx}}{v_x}
$$
est la somme d'une fonction de $x$ et d'une fonction de $y$. Donc $h_{xy}=0$, 
ce qui donne :
$$
v_xv_{xx}\left(\fr{w_{yyy}}{w_y}\right)_{\!\!y}  =  w_yw_{yy} \left(\fr{v_{xxx}}{v_x}\right)_{\!\!x},
$$
ou mieux :
$$
\fr{1}{w_yw_{yy}}\left(\fr{w_{yyy}}{w_y}\right)_{\!\!y}  
= 
\fr{1}{v_xv_{xx}} \left(\fr{v_{xxx}}{v_x}\right)_x
\equiv 
4p\in \C,
$$
puisque le premier membre ne d\'epend que de $y$ et le deuxi\`eme que de $x$.

\sk
On a obtenu que $v_x$ et $w_y$, comme fonctions d'une variable, sont solutions de 
la m\^eme \'equation diff\'erentielle :
\beq
\label{edo}
\left( \fr{z''}{z} \right)' = 4 pzz'.
\eeq
Si $z\not\equiv 0$ v\'erifie (\ref{edo}), $z'' = 2pz^3 + qz$ pour un 
$q\in \C$, donc $z''z' = 2pz^3z' + qzz'$ et finalement :
\beq
\label{edo'}
(z')^2  = pz^4  + qz^2 + r, \qquad q,r\in \C.
\eeq
R\'eciproquement, supposons que  :
$$
v_{xx}^2  = pv_x^4  +  qv_x^2 + r  \, ; \;\; w_{yy}^2  = pw_y^4  +  q'w_y^2 + r',
$$
donc :
$$
\fr{v_{xxx}}{v_x} = 2pv_x^2 + q  \, ;  \;\; \fr{w_{yyy}}{w_y} = 2pw_y^2  + q'.
$$
L'\'equation (\ref{eq1ter}) est v\'erifi\'ee si et seulement si :
$$
( (2pv_x^2 + q) + (2pw_y^2 + q') )(v_x^2 - w_y^2) = 2( (pv_x^4 + qv_x^2 + r) - (pw_y^4 + q'w_y^2 + r') ).
$$
Apr\`es simplification, on obtient $(q-q')v_x^2 + 2r = (q'-q)w_y^2 + 2r'$ et donc 
$q=q'$, $r=r'$. 

Finalement, si $q=q'$ et $r=r'$, le premier membre de l'\'equation (\ref{eq2}) devient :
$$
X_u \left( \fr{ v_{xxx}/v_x - w_{yyy}/w_y }{v_x^2-w_y^2} \right) = 2pX_u1 = 0.
$$
L'\'equation (\ref{eq2}) est v\'erifi\'ee. En r\'esum\'e, on a d\'emontr\'e :
\bpr
\label{les4cas}
Soit $u(x,y) = v(x) + w(y)$ un germe de fonction  qui 
v\'erifie (\ref{transversebis}). Le tissu $\cal{T}[u]$ est 
de rang maximal si et seulement si le tissu 
$\cal{T}(\, x+y, \, x-y, \, u(x,y)\,)$ est de rang maximal, si et seulement si 
l'une des conditions suivantes est v\'erifi\'ee :
\be
\item $v_x$ et $w_y$ sont constantes ;
\item $v_x$ est une constante $c\in \C^*$ et $w_{yy} = k(c^2 - w_y^2)$ pour un $k\in \C^*$ ;
\item $w_y$ est une constante $c\in \C^*$ et $v_{xx} = k(c^2 - v_x^2)$ pour un $k\in \C^*$ ;
\item il existe $p,q,r\in \C$ (non tous nuls) tels que $v_x$ et $w_y$ sont des solutions 
non constantes de l'\'equation diff\'erentielle $(z')^2 = pz^4 + qz^2 + r$.
\ee
\epr

\bpf[D\'emonstration du Th\'eor\`eme \ref{tous}]
C'est un exercice assez r\'ebarbatif qui consiste \`a 
r\'e\'ecrire les \'equations de la proposition, 
en distinguant une petite multitude de cas, et \`a 
v\'erifier qu'\`a chaque fois la solution $u(x,y)$ d\'efinit,
\`a l'action pr\`es du groupe ${\rm Sym}\,(\cal{T}_0)$,
un tissu du Th\'eor\`eme \ref{famille} ou de l'un des 
types (A), (B), (C), (D) ou (E) du Th\'eor\`eme \ref{autres}.
Les fonctions $u$ qui vont intervenir sont des fonctions 
classiques. On ne fera pas r\'ef\'erence au point de base,
qui peut \^etre n'importe quel point tel que 
(\ref{transversebis}) soit v\'erifi\'e.

\sk
Soit $u(x,y)$ et $u'(x,y)$ deux fonctions du type consid\'er\'e 
dans la Proposition \ref{les4cas}. On suppose de plus 
que les feuilletages qu'elles d\'efinissent ne sont pas des
faisceaux de droites. Sous cette condition, on sait qu'elles d\'efinissent 
le m\^eme tissu si et seulement si elles d\'efinissent le m\^eme feuilletage.
C'est le cas si et seulement si $u'$ est de la forme 
$$
u'(x,y) = cu(x,y) + e, \qquad (c\in \C^*, \; e\in \C).
$$
(En effet, la forme impos\'ee $u(x,y) = v(x) + w(y)$ est une relation 
ab\'elienne non triviale du $3$-tissu $\cal{T}(\,x,\,y,\, u(x,y))$. Toute 
autre relation de ce tissu est une combinaison lin\'eaire de 
celle-ci et d'une \og relation constante \fg.)

\sk
On utilisera la notation (\ref{dilat}).

\sk
\nk
{\em 1. --- } Si $u$ v\'erifie la premi\`ere propri\'et\'e de la proposition, 
$u$ est une fonction affine, $\cal{F}(u)$ un faisceau de droites parall\`eles
et $\cal{T}[u]$ un tissu alg\'ebrique.

\sk
\nk
{\em 2. --- }  Si $u$ v\'erifie la deuxi\`eme propri\'et\'e de la proposition,
quitte \`a remplacer $u$ par $u' \sim u$, 
on se ram\`ene au cas o\`u $c=1$, $k=1$. 
L'\'equation $w_{yy} = 1 - w_y^2$ a la solution $w_y = \tanh y$. 
On obtient $u(x,y) \sim x + \log \cosh y$ et des tissus du type (C).

\sk
\nk
{\em 3. --- }  Si $u$ v\'erifie la troisi\`eme propri\'et\'e de la proposition,
la sym\'etrie $(x,y)\mapsto (y,x)$ permet de se ramener au cas pr\'ec\'edent.
 
\sk
\nk
{\em 4. --- } On suppose maintenant que $u$ v\'erifie la quatri\`eme propri\'et\'e. On note :
$$
{\rm E(p,q,r)} : \qquad (z')^2 = pz^4 + qz^2 + r,
$$
l'\'equation diff\'erentielle v\'erifi\'ee par $v_x = u_x$ et $w_y = u_y$.
Si $z$ est une solution particuli\`ere de l'\'equation pr\'ec\'edente,
la solution g\'en\'erale est $t \mapsto \pm z(t+t')$. Donc, \'etant donn\'e 
{\em une} primitive $v$ de $z$, on a :
$$
u(x,y) \sim  v(x) \pm v(y).
$$
Selon les cas les tissus associ\'es \`a des signes diff\'erents sont \'equivalents
ou ne le sont pas.

\sk
D'autre part, en rempla\c{c}ant $u$ par $u'\sim u$, on peut remplacer l'\'equation 
E(p,q,r) par l'\'equation E($c^2d^2$p, $c^2$q$/d^2$,r$/(c^2d^2)$)
pour tout $c,d\in \C^*$.

\sk
\nk
{\em 4-a. --- } Si $p=0$ et $q=0$, on se ram\`ene \`a $r=2$, d'o\`u $z(t)=t/2$ et $v(x)=x^2$. On obtient 
$u(x,y) \sim x^2 \pm y^2$ et, selon le signe, des tissus du type (D) ou du type (E).

\sk
\nk
{\em 4-b. --- } Si $p=0$ et $r=0$, on se ram\`ene \`a $q=1$, d'o\`u $z(t) = \e{t}$, $v(x) = \e{x}$. On obtient 
$u(x,y) \sim \e{x} \pm \e{y}$. Une translation permet de passer du signe $-$ au signe $+$.
On obtient des tissus du type (C).

\sk
\nk
{\em 4-c. --- } Si $p=0$ et  $qr\neq 0$, on se ram\`ene \`a $q=1$, $r=1$, $z(t)=\sinh t$, $v(x)=\cosh x$.
On obtient $u(x,y) \sim \cosh x \pm \cosh y$. 
Une translation permet de passer du signe $-$ au signe $+$.
On obtient des tissus du type (B).

\sk
\nk
{\em 4-d. --- } Si $q=0$ et $r=0$, on se ram\`ene \`a $p=1$, d'o\`u $z(t)= 1/t$, $v(x) = \log x$. On obtient :
$u(x,y) \sim \log x \pm \log y$. Si le signe est $-$, on obtient un faisceau de droites,
donc un tissu alg\'ebrique. 
Si le signe est $+$, on obtient des tissus du type (D).

\sk
\nk
{\em 4-e. --- } On suppose maintenant $p\neq 0$, $(q,r)\neq (0,0)$.  En rempla\c{c}ant 
$u$ par $u'\sim u$, on se ram\`ene au cas o\`u l'\'equation 
E(p,q,r) est de la forme suivante :
\beq
\label{Ek}
{\rm E(k)} : \qquad (z')^2 = (1-z^2)(1-k^2z^2).
\eeq
Cette \'equation a comme solution particuli\`ere 
la fonction de Jacobi ${\rm sn}_k t$
si $k\neq 0,1$, la fonction $\sin t$ si $k=0$ et la fonction $\tanh t$
si $k=1$. On a 
$$
\int {\rm sn}_k t \, dt = \fr{1}{k} \log ({\rm dn}_k t - {\rm cn}_k t) + {\rm Cte} \sim \log ({\rm sn_l} t), 
$$
o\`u $l=(1-k)/(1+k)$ (voir \cite{La} (2.7.1)).

\sk
\nk
{\em 4-e (i). --- } Si $k^2=0$, on obtient $u(x,y) \sim \log \sin x \pm \log \sin y$.
Selon que le signe est $+$ ou $-$, les tissus obtenus sont du type (B) 
ou du type (A).

\sk
\nk
{\em 4-e (j). --- } Si $k^2=1$, on obtient $u(x,y) \sim \log \tanh x \pm \log \tanh y$.
Une translation permet de passer du signe $-$ au signe $+$.
On obtient des tissus du type (A).

\sk
\nk
{\em 4-e (k). --- } Si $k^2\neq 0,1$, on obtient $u(x,y) \sim \log \sn_l x \pm \log \sn_l y$.
Compte tenu de (\ref{inverse-sn}), une translation permet de passer du signe $-$ au signe $+$.
On obtient les tissus du Th\'eor\`eme \ref{famille}.

\epf

\end{document}